\documentclass[10pt, a4paper, reqno, english]{amsart}

\usepackage{amssymb, mathrsfs, mathtools}

\usepackage[T1]{fontenc}
\usepackage[hidelinks]{hyperref}
\usepackage{microtype}
\usepackage{tikz}
\usetikzlibrary{shapes.geometric}
\usepackage{tikz-cd}
\usepackage{lmodern}
\usepackage{todonotes}

\theoremstyle{plain}
\newtheorem*{theorem*}{Theorem}
\newtheorem{prop}{Proposition}
\newtheorem{lemma}[prop]{Lemma}
\newtheorem{theorem}[prop]{Theorem}
\theoremstyle{remark}

\numberwithin{prop}{section}
\numberwithin{equation}{section}

\newcommand\ab{\boldsymbol{a}}
\newcommand\Bb{\boldsymbol{B}}
\newcommand\Kb{\boldsymbol{K}}
\newcommand\oneb{\boldsymbol{1}}
\newcommand\ddd{\,\mathrm{d}}
\newcommand\PP{\mathbb{P}}
\newcommand\QQ{\mathbb{Q}}
\newcommand\RR{\mathbb{R}}
\newcommand\AAA{\mathbb{A}}
\newcommand\CC{\mathbb{C}}
\newcommand\ZZ{\mathbb{Z}}
\newcommand{\GG}{\mathbb{G}}
\newcommand{\GGm}{\GG_{\mathrm{m}}}
\newcommand{\GGmOK}{\GG_{\mathrm{m},\OK}}
\newcommand\afr{{\mathfrak{a}}}
\newcommand\cfr{{\mathfrak{c}}}
\newcommand\pfr{{\mathfrak{p}}}
\newcommand\qfr{{\mathfrak{q}}}
\newcommand\N{\mathfrak{N}}
\newcommand\afrb{{\underline{\mathfrak{a}}}}
\newcommand\cfrb{{\underline{\mathfrak{c}}}}
\newcommand{\Cs}{\mathscr{C}}
\newcommand{\Is}{\mathscr{I}}
\newcommand{\IK}{\Is_K}
\newcommand{\Fs}{\mathscr{F}}
\newcommand{\Gs}{\mathscr{G}}
\newcommand{\Cmax}{\Cs^{\mathrm{max}}}
\newcommand{\Os}{\mathcal{O}}
\newcommand{\OK}{\Os_K}
\newcommand{\Op}{\Os_{\pfr}}
\newcommand\Ac{{\mathcal{A}}}
\newcommand\Dc{{\mathcal{D}}}
\newcommand\Sc{{\mathcal{S}}}
\newcommand\Sct{{\widetilde{\mathcal{S}}}}
\newcommand\Sp{{\widetilde{\mathcal{S}}_\pfr}}
\newcommand\Uc{{\mathcal{U}}}
\newcommand\Uct{{\widetilde{\mathcal{U}}}}
\newcommand\Yc{{\mathcal{Y}}}

\newcommand\cYc{{\prescript{}{\cfrb}{\mathcal{Y}}}}
\newcommand\cYcU{\prescript{}{\cfrb}{\mathcal{Y}}_{\Uct}}
\newcommand\Tc{{\mathcal{T}}}
\newcommand\Rc{{\mathcal{R}}}
\newcommand\Lc{{\mathcal{L}}}

\newcommand\crho{{\prescript{}{\cfrb}{\rho}}}
\newcommand\Nt{\widetilde{N}}
\newcommand\St{{\widetilde{S}}}
\newcommand\Vt{\widetilde{V}}
\newcommand\sums[1]{\sum_{\substack{#1}}}
\newcommand\ints[1]{\int_{\substack{#1}}}
\newcommand\bigwhere[2]{\left\{#1:\ \begin{aligned}#2\end{aligned}\right\}}
\newcommand\congr[3]{#1 \equiv #2 \pmod{#3}}
\newcommand\fin{\mathrm{fin}}

\DeclareMathOperator{\vol}{vol}
\DeclareMathOperator{\Pic}{Pic}
\DeclareMathOperator{\Spec}{Spec}
\DeclareMathOperator{\Cl}{Cl}
\DeclareMathOperator{\Eff}{Eff}
\DeclareMathOperator{\rk}{rk}

\AtBeginDocument{%
   \def\MR#1{}
}

\begin{document}

\title[Integral points over number fields]
{Integral points over number fields:\\a Clemens complex jigsaw puzzle}

\author[C.~Bernert]{Christian Bernert}

\address{Institute of Science and Technology Austria, Am Campus 1, 3400 Klosterneuburg, Austria}

\email{christian.bernert@ist.ac.at}

\author[U.~Derenthal]{Ulrich Derenthal} 

\address{Institut f\"ur Algebra, Zahlentheorie und Diskrete Mathematik, Leibniz Universit\"at Hannover, Welfengarten 1, 30167 Hannover, Germany}

\email{\href{mailto:derenthal@math.uni-hannover.de}{derenthal@math.uni-hannover.de}}

\author[J.~Ortmann]{Judith Ortmann}

\address{Institut f\"ur Algebra, Zahlentheorie und Diskrete Mathematik, Leibniz Universit\"at Hannover, Welfengarten 1, 30167 Hannover, Germany}

\email{\href{mailto:ortmann@math.uni-hannover.de}{ortmann@math.uni-hannover.de}}

\author[F.~Wilsch]{Florian Wilsch}
\address{Mathematisches Institut, Georg-August-Universität Göttingen, Bunsenstra\ss{}e 3--5, 37073 Göttingen, Germany}
\email{\href{mailto:florian.wilsch@mathematik.uni-goettingen.de}{florian.wilsch@mathematik.uni-goettingen.de}}

\address{}

\date{September 9, 2026}

%11G35: NT -> Arithmetic algebraic geometry (Diophantine geometry) -> Varieties over global fields
%11D45: NT -> Diophantine equations -> Counting solutions of Diophantine equations
%14G05: AG -> Arithmetic problems in algebraic geometry; Diophantine geometry -> Rational points
%52B20: Convex and discrete geometry -> Polytopes and polyhedra -> Lattice polytopes in convex geometry (including relations with commutative algebra and algebraic geometry) [See also 06A11, 13F20, 13F55, 13Hxx, 52C05, 52C07]

\setcounter{tocdepth}{1}

\begin{abstract}
  We prove an asymptotic formula for the number of integral points of bounded log anticanonical height on a singular quartic del Pezzo surface over arbitrary number fields, with respect to the largest admissible boundary divisor. The resulting Clemens complex is more complicated than usual, and leads to particularly interesting effective cone constants, associated with exponentially many polytopes whose volumes appear in the expected formula. Like a jigsaw puzzle, these polytopes fit together to one large polytope. The volume of this polytope appears in the asymptotic formula that we obtain using the universal torsor method via o-minimal structures.

\end{abstract}
  
\maketitle
\tableofcontents

\section{Introduction}

The Manin--Peyre conjecture predicts the asymptotic behavior of the number of rational points of bounded height on a Fano variety. While a proof in the general case remains elusive, our understanding of the constants appearing in this asymptotic formula has consolidated over the past decades, in particular conforming to both the results obtained by harmonic analysis techniques as well as by the universal torsor method.

On the other hand, the analogous situation for integral points is much less streamlined. Only recently \cite{Wil24,Santens23} (based on \cite{CLT10,CLT12,CLTtoric}), a precise description of the expected main term has emerged. As integral points on a log Fano variety with respect to a boundary divisor $D$ tend to accumulate near the minimal strata of $D$, it transpires that the  distribution of integral points is governed by the maximal faces of the associated analytic Clemens complex.
From the point of view of harmonic analysis, as worked out in \cite{Santens23} in the case of toric varieties, it is then natural to study the contribution from each face separately, thus resulting in an expected formula that involves a sum over these maximal faces.

In an application of the torsor method, however, one tends to obtain an asymptotic formula that naturally involves only one main term. As both the predicted contributions from each face of the Clemens complex as well as the constant obtained from the torsor method involve volumes of certain polytopes, it is natural to wonder to which extent these pictures can be unified.

This article explores this phenomenon by revisiting the singular quartic del Pezzo surface with an $A_1$ and $A_3$ singularity, which has served as an instructive example in other contexts before: for rational points over $\QQ$ \cite{D09}, imaginary quadratic fields \cite{DF14b}, arbitrary number fields \cite{FP16,DP20} and function fields \cite{Bourqui}, and for integral points over $\QQ$ \cite{DW24} and imaginary quadratic fields \cite[Theorem~2.7 (cases 1, 4)]{OrtmannThesis}.
It is neither toric \cite[Remark~6]{D14} nor an equivariant compactification of a vector group \cite{DL10}.

\begin{figure}[ht]
\begin{tikzpicture}
  \def\pcolora{gray!90}
  \def\pcolorb{gray!20}
  \def\popacity{0.4}

  \path[fill=green, opacity=0.5] (1, 0) -- (1, 5) -- (0, 2) -- (0, 0);
  \path[fill=blue,  opacity=0.5] (0, 2) -- (1, 5) -- (0, 3);
  \path[fill=red,   opacity=0.5] (0, 3) -- (1, 5) -- (0, 4);
  \path[fill=orange,opacity=0.5] (0, 4) -- (1, 5) -- (0, 5);

  \path[fill=\pcolora, opacity=0.0] (0, 0) -- (0, 5) -- (1, 5) -- (1,3);
  \path[fill=\pcolora, opacity=0.25] (0, 0) -- (1, 3) -- (1, 2);
  \path[fill=\pcolora, opacity=0.5] (0, 0) -- (1, 2) -- (1, 1);
  \path[fill=\pcolora, opacity=0.75] (0, 0) -- (1, 1) -- (1, 0);

  % for B_1
  \draw[thick] (0, 0) -- (0, 5);
  \draw[thick] (0, 0) -- (1, 3);
  \draw[thick] (0, 0) -- (1, 2);
  \draw[thick] (0, 0) -- (1, 1);
  \draw[thick] (0, 0) -- (1, 0);
  % for B_2
  \draw[thick] (1, 0) -- (1, 5);
  \draw[thick] (0, 2) -- (1, 5);
  \draw[thick] (0, 3) -- (1, 5);
  \draw[thick] (0, 4) -- (1, 5);
  \draw[thick] (0, 5) -- (1, 5);
 \end{tikzpicture}
 \begin{tikzpicture}
  \def\pcolora{gray!90}
  \def\pcolorb{gray!20}
  \def\popacity{0.4}

  \path[fill=green, opacity=0.5] (2, 0) -- (2, 5) -- (1/3, 0);
  \path[fill=blue,  opacity=0.5] (1/3,0) --(2, 5) -- (0, 1) -- (0, 0);
  \path[fill=red,   opacity=0.5] (0, 1) -- (2, 5) -- (0, 3);
  \path[fill=orange,opacity=0.5] (0, 3) -- (2, 5) -- (0, 5);

  \path[fill=\pcolora, opacity=0.0] (0, 0) -- (0, 5) -- (5/3, 5);
  \path[fill=\pcolora, opacity=0.25] (0, 0) -- (5/3, 5) -- (2,5) -- (2, 4);
  \path[fill=\pcolora, opacity=0.5] (0, 0) -- (2, 4) -- (2, 2);
  \path[fill=\pcolora, opacity=0.75] (0, 0) -- (2, 2) -- (2, 0);

  % for B_1
  \draw[thick] (0, 0) -- (0, 5);
  \draw[thick] (0, 0) -- (5/3, 5);
  \draw[thick] (0, 0) -- (2, 4);
  \draw[thick] (0, 0) -- (2, 2);
  \draw[thick] (0, 0) -- (2, 0);
  % for B_2
  \draw[thick] (2,   0) -- (2, 5);
  \draw[thick] (1/3, 0) -- (2, 5);
  \draw[thick] (0,   1) -- (2, 5);
  \draw[thick] (0,   3) -- (2, 5);
  \draw[thick] (0,   5) -- (2, 5);
 \end{tikzpicture}
 \begin{tikzpicture}
  \def\pcolora{gray!90}
  \def\pcolorb{gray!20}
  \def\popacity{0.4}
  \path[fill=green, opacity=0.5] (3, 0) -- (3, 5) -- (4/3, 0);
  \path[fill=blue,  opacity=0.5] (4/3,0) --(3, 5) -- (1/2, 0);
  \path[fill=red,   opacity=0.5] (1/2, 0) -- (3, 5) -- (0, 2) -- (0, 0);
  \path[fill=orange,opacity=0.5] (0, 2) -- (3, 5) -- (0, 5);

  \path[fill=\pcolora, opacity= 0.0] (0, 0) -- (0, 5) -- (5/3, 5);
  \path[fill=\pcolora, opacity=0.25] (0, 0) -- (5/3, 5) -- (5/2, 5);
  \path[fill=\pcolora, opacity= 0.5] (0, 0) -- (5/2, 5) -- (3, 5) -- (3, 3);
  \path[fill=\pcolora, opacity=0.75] (0, 0) -- (3, 3) -- (3, 0);
  % for B_1
  \draw[thick] (0, 0) -- (0, 5);
  \draw[thick] (0, 0) -- (5/3, 5);
  \draw[thick] (0, 0) -- (5/2, 5);
  \draw[thick] (0, 0) -- (3, 3);
  \draw[thick] (0, 0) -- (3, 0);
  % for B_2
  \draw[thick] (3,   0) -- (3, 5);
  \draw[thick] (4/3, 0) -- (3, 5);
  \draw[thick] (1/2,   0) -- (3, 5);
  \draw[thick] (0,   2) -- (3, 5);
  \draw[thick] (0,   5) -- (3, 5);
 \end{tikzpicture}
 
 \caption{The jigsaw puzzle for a number field with two archimedean places.}\label{fig:jigsaw}
 
\end{figure}
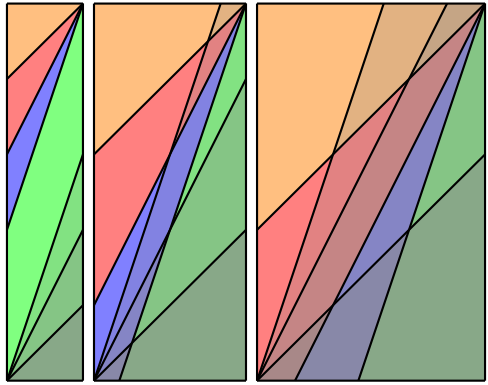

Before stating our results, let us explain why this ``jigsaw puzzle'' phenomenon has not been a prominent feature in previous work. Indeed, the complexity of the Clemens complex is governed by two aspects, the choice of the boundary divisor and of the ground field:

If the boundary $D$ is a prime divisor (as in \cite{BD25integral}) or, more generally, if all divisor components intersect over the base field (as for wonderful compactifications~\cite{TBT,Chow}), the analytic Clemens complex has only one maximal face, even over an arbitrary number field $K$. Meanwhile, when $D$ is arbitrary and one works over $\QQ$ (as in \cite{DW24}) or an imaginary quadratic field (as in \cite{OrtmannPaper}), the number of maximal faces tends to be at most the number of irreducible components of $D$, especially for the well-studied case of del Pezzo surfaces (though it can be much larger in general as the possible ways to intersect divisors grows with the dimension).

For our example, we have therefore chosen the largest admissible \cite[Theorem~10]{DW24} boundary divisor $D$, with five components and four maximal faces, and we work over an arbitrary number field. In this case, the analytic Clemens complex has $4^{q+1}$ maximal faces, where $q$ is the rank of the unit group of $K$ (so that $q+1$ is its number of archimedean places). Figure~\ref{fig:jigsaw} shows three cross sections of the four-dimensional polytope in the case $q=1$; see Section~\ref{sec:alpha} for further details.

\subsection{The counting result}

We now set out to describe our results in more detail. Let $K$ be a number field. Consider the singular quartic del Pezzo surface $S \subset \PP^4_K$ defined by
\begin{equation}\label{eq:surface}
    x_0x_3-x_2x_4 = x_0x_1 + x_1x_3 + x_2^2 = 0.
\end{equation}
It admits two singular points: $Q_1=(0:1:0:0:0)$ as an $A_1$-singularity and $Q_2 = (0:0:0:0:1)$ as an $A_3$-singularity.
We are interested in integral points with respect to the line $L = \{x_0=x_2=x_3=0\}$ through the two singularities (which corresponds to a divisor $D$ on the minimal desingularization of $S$ that has the maximal number of five components among those that define a \emph{weak del Pezzo pair}~\cite[Thm.~10, Lem.~12]{DW24}).
The surface $S$ contains two more lines: $L'=\{x_1=x_2=x_3=0\}$ and $L''=\{x_0=x_1=x_2=0\}$, both passing through $Q_2$; let $V \subset S$ be the complement of the three lines.
Let $\Sc \subset \PP^4_{\OK}$ be the integral model of $S$ defined by the same equations~\eqref{eq:surface} over $\OK$, let $\Uc = \Sc \setminus \overline{L}$ be the complement of the Zariski closure $\overline{L}$ of $L$ in $\Sc$, and let $U = S \setminus L$ be its generic fiber.

We consider the following height function (which turns out to correspond to a log anticanonical height on the desingularization of $S$, see Section~\ref{sec:height}): For $x = (x_0:\dots:x_4) \in S(K)$ outside the lines, put
\begin{equation}\label{eq:height-intro}
    H(x)=\prod_{v \in \Omega_K} \max\{|x_0|_v,|x_2|_v,|x_3|_v\} = \frac{\prod_{v \mid \infty} \max\{|x_0|_v,|x_2|_v,|x_3|_v\}}{\N(x_0\OK+x_2\OK+x_3\OK)}.
\end{equation}

Our first main result determines the asymptotic behavior of integral points with respect to the boundary $L$ of bounded height outside the lines:
\begin{equation*}
    N_{\Uc,V,H}(B) \coloneq |\{x \in \Uc(\OK) \cap V(K) : H(x) \le B\}|.
\end{equation*}

\begin{theorem}\label{thm:main_concrete}
    We have
    \begin{equation*}
        N_{\Uc,V,H}(B) = c B(\log B)^{2+2q} + O(B(\log B)^{1+2q}\log\log B),
    \end{equation*}
    with
    \begin{equation}\label{eq:c}
        c = \alpha \frac{\rho_K}{|\Delta_K|} \prod_{v \in \Omega_K} \omega_v,
    \end{equation}
    where $\rho_K$ as in \eqref{eq:def_rho_K} is the residue of the Dedekind zeta function $\zeta_K$ at $s=1$, and 
    \begin{align}
        \label{eq:alpha}\alpha &= \frac{1}{q!^2} \int_{\substack{t_1,t_2 \ge 0\\t_1+t_2 \le 1}} (1-t_1-t_2)^q \ddd t_1 \ddd t_2 = \frac{1}{q!(q+2)!},\\
        \label{eq:omega_v}\omega_v &=
        \begin{cases}
            4, &\text{if $v$ is real,}\\
            4\pi^2, &\text{if $v$ is complex,}\\
            1-\frac 1{\N\pfr^2}, &\text{if $v = \pfr$ is finite}.
        \end{cases}
    \end{align}
\end{theorem}

The counting argument ends up being fairly short, even though we work over arbitrary number fields. This is partly due to the way the variables are restricted in Section~\ref{sec:restrictions}, a new ingredient to the proof that simplifies the analysis in Section~\ref{sec:main_contribution} considerably and would have a similar effect on the analytic arguments in previous works on this surface (such as \cite[\S\S\,7, 10, 11]{FP16}, for example). It is also convenient that we do not need a Möbius inversion.

Our second goal in this paper is to show that this result agrees with the asymptotic formula predicted by \cite{Wil24,Santens23}. To describe it, we first need to introduce some notation.

\subsection{The expected formula}

As usual, the asymptotic formula should be interpreted after resolving the singularities. To this end, let $\pi \colon \St \to S$ be the minimal desingularization. Let $\Vt = \pi^{-1}(V)$; this is the complement of the negative curves on $\St$. As $U$ lies within the smooth locus, it remains unchanged under $\pi$, and we liberally identify it with its preimage. Write $D = \St \setminus U$ for the boundary divisor; it has five components (the strict transform of the line and $3+1$ exceptional curves lying above the singular points) that we will describe in greater detail in Section~\ref{sec:torsor}.

\begin{figure}[ht]
  \begin{center}
    \begin{tikzpicture}
      \node[label=$A_7$] (A7) at (  -3, 0) [draw, shape=circle, fill=black, scale=.3]{};
      \node[label=$A_5$] (A5) at (-1.5, 0) [draw, shape=circle, fill=black, scale=.3]{};
      \node[label=$A_4$] (A4) at (   0, 0) [draw, shape=circle, fill=black, scale=.3]{};
      \node[label=$A_3$] (A3) at ( 1.5, 0) [draw, shape=circle, fill=black, scale=.3]{};
      \node[label=$A_6$] (A6) at (   3, 0) [draw, shape=circle, fill=black, scale=.3]{};

      \draw[thick] (A7) -- (A5) node[pos=0.5, anchor=north]{$(57)$};
      \draw[thick] (A5) -- (A4) node[pos=0.5, anchor=north]{$(45)$};
      \draw[thick] (A4) -- (A3) node[pos=0.5, anchor=north]{$(34)$};
      \draw[thick] (A3) -- (A6) node[pos=0.5, anchor=north]{$(36)$};
    \end{tikzpicture}
    \caption{The $K_v$-analytic Clemens complex of $D$ is a path on five vertices and does not depend on $v$. Each vertex corresponds to an irreducible component of $D$ on a minimal desingularization, labeled as in Section~\ref{sec:torsor}.}\label{fig:clemens_complex}
  \end{center}
\end{figure}
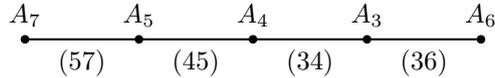

The expected asymptotic formula heavily features the \emph{Clemens complex}, which encodes incidence properties of the divisor $D$ over each completion of the base field $K$. Owing to the splitness of $\St$, the \emph{$K_v$-analytic Clemens complex}~\cite[\S\,3.1]{CLT10} $\Cs_v(D)$ is in fact independent of the place $v$ appearing in its definition. It consists of one vertex for each of the five irreducible components of $D$ with an edge added between the vertices associated with intersecting components; in other words, it is the full subgraph of the Dynkin diagram on the components of $D$. The complex, including a labeling of its faces, is given in Figure~\ref{fig:clemens_complex}: each vertex is named after the corresponding divisor component, and each one-dimensional (maximal) face is labeled by the tuple of the indices corresponding to its two vertices.
The \emph{analytic Clemens complex} is the product $\Cs (D) = \prod_{v\mid\infty} \Cs_v(D)$. These Clemens complexes are partially ordered via incidence; write $\prec$ for this order. Maximal faces $\Bb\in \Cmax(D)$ are tuples $(B_v)_{v\mid\infty}$ of maximal faces $B_v\in \Cmax_v(D)$.

Concretely, the expected formula 
\begin{equation}\label{eq:expected-formula}
        N_{\Uc,V,H}(B) = c_\infty c_\fin B(\log B)^{b-1} (1+o(1)),
    \end{equation}
described by the fourth author~\cite[\S\,2.5]{Wil24} involves two factors
\begin{align*}
    c_\infty &=
    \frac{1}{|\Delta_K|} \sum_{\Bb\in \Cmax(D)}\alpha_{\Bb}(\St) \prod_{v \mid \infty} \tau_{{B_v},v}(Z_{B_v}(K_v)),\\
    c_\fin &=
    \rho_K^{\rk \Pic U} \prod_{\pfr} \left(1-\frac{1}{\N\pfr}\right)^{\rk \Pic U} \tau_{(\St,D),\pfr}(\Uc(\Os_\pfr))
\end{align*}
of the leading constant and $b = \rk\Pic U +\sum_{v\mid\infty}(\dim \Cs_v(D)+1)$ in the exponent of $\log B$. Santens~\cite[Theorem~6.11, \S\,6.3]{Santens23} shows that this expectation agrees with his more general conjecture.

We will explicitly compute these ingredients of the expected formula in Section~\ref{sec:expected}. In particular, this includes solving our jigsaw puzzle.

\begin{prop}\label{prop:puzzle}
    With $\alpha$ as in \eqref{eq:alpha}, we have
    \begin{equation*}
        \alpha = \sum_{\Bb\in \Cmax(D)}\alpha_{\Bb}(\St).
    \end{equation*}
\end{prop}

We thus obtain the following reformulation of Theorem~\ref{thm:main_concrete}.

\begin{theorem}\label{thm:counting-expected-formula}
    The formula~\eqref{eq:expected-formula} holds true, confirming both expectations~\cite{Wil24,Santens23}.
\end{theorem}

\subsection{Overview of the proof}

In Section \ref{sec:torsor}, we lift the set of integral points to twisted universal torsors. This is the igniting spark for both the parameterization underlying the counting techniques used in the proof of Theorem \ref{thm:main_concrete} and the computation of the Tamagawa measures in the expected formula, as required for the proof of Theorem \ref{thm:counting-expected-formula}.
In Section~\ref{sec:expected}, we deduce Theorem~\ref{thm:counting-expected-formula} from Theorem \ref{thm:main_concrete}. As explained in the introduction, a key novelty here is a comparison of the effective cone constants in Proposition~\ref{prop:puzzle}.
In Section~\ref{sec:parameterization}, we initiate the proof of Theorem~\ref{thm:main_concrete} by setting up a concrete counting problem on the universal torsors. Working over number fields, this requires the construction of a suitable fundamental domain. In Proposition~\ref{prop:restrict_a6_a7}, we also prepare for the counting by introducing a truncation of some of the torsor variables. This move, motivated by \cite{BD24}, turns out to be very useful for estimating the error terms in the later stages of the counting argument.
Finally, in Section \ref{sec:main_contribution}, we reduce to a lattice point counting problem which we can solve using o-minimal structures. After computing the local densities and summing over the remaining variables, this eventually leads to a proof of Theorem \ref{thm:main_concrete}.

This work is based on the third author's Ph.D. thesis \cite[Theorem~2.7, case~3]{OrtmannThesis}, which establishes the asymptotic formula in Theorem~\ref{thm:main_concrete}. New ingredients in our present work include the interpretation of the constant $\alpha$ as in Theorem~\ref{thm:counting-expected-formula}, the more conceptual computation of the nonarchimedean densities using model measures and numbers of points modulo $\pfr$ (instead of $\pfr$-adic integrals) in Section~\ref{sec:tamagawa}, as well as a simplified counting machinery via the restrictions in Section~\ref{sec:restrictions}.

\subsection{Notation and conventions}

We use the standard notation for number fields as in \cite[\S\,1.5]{BD24}, in particular its ring of integers $\OK$, its regulator $R_K$, its class number $h_K$, its discriminant $\Delta_K$ and its number of roots of unity $|\mu_K|$. Let $r_1$ be the number of its real embeddings, and $r_2$ the number of pairs of complex embeddings. Then $q\coloneq r_1+r_2-1$ is the rank of $\OK^\times$. Let
\begin{equation}\label{eq:def_rho_K}
    \rho_K\coloneq \frac{2^{r_1}(2\pi)^{r_2}R_K h_K}{|\mu_K|\cdot|\Delta_K|^{1/2}}.
\end{equation}
Let $d = [K:\QQ]$, $\IK$ be the monoid of nonzero ideals in $\OK$, $U_K$ the subgroup of $\OK^\times$ generated by a chosen system of fundamental units, and $\Cl_K$ the ideal class group.

Let $\Omega_K$ be the set of places of $K$. For $v$ lying in its subset $\Omega_\infty$ of archimedean places, we also write $v \mid \infty$. For $v \in \Omega_K$ above $w \in \Omega_\QQ$, let $d_v = [K_v : \QQ_w]$, and $|\cdot|_v = |N_{K_v/\QQ_w}(\cdot)|_w$ (with $|\cdot|_w$ the usual $p$-adic or real absolute value on $\QQ_w$). For $a \in K$, let $a^{(v)}$ be its image in the completion $K_v$, let $|a|_v = |a^{(v)}|_v$, and let $\sigma(a) = (a^{(v)})_{v \mid \infty} \in \prod_{v \mid \infty} K_v$.

Let $\N\qfr$ be the norm of a fractional ideal $\qfr$. The letter $\pfr$ always denotes a nonzero prime ideal in $\OK$ corresponding to a finite place $v$, and products over $\pfr$ run through all of them (possibly subject to further conditions as stated). Let $\Op$ be the valuation ring in $K_\pfr=K_v$.

When we use Vinogradov's $\ll$-notation or Landau's $O$-notation, the corresponding inequalities are meant to hold for all values in the relevant range, and the implied constants may depend only on $K$. We write $X_1 \asymp X_2$ for $X_1 \ll X_2 \ll X_1$.
Volumes of subsets of $\RR^n$ or $\CC^n \cong \RR^{2n}$ are computed with respect to the usual Lebesgue measure, unless stated otherwise.
For our height bound $B$, we assume $B \ge 3$.

\subsection*{Acknowledgements}
We would like to thank the referee for their careful reading and helpful remarks that helped to improve this article.
The first two authors were supported by the Deutsche For\-schungs\-gemeinschaft (DFG) -- 512730679 (RTG2965). The third author was partially supported by the Caroline Herschel Programme of Leibniz University Hannover.
The fourth author was supported by the Deutsche Forschungsgemeinschaft (DFG) -- 398436923 (RTG2491).

\section{Passage to universal torsors}\label{sec:torsor}

\subsection{Universal torsors over $K$} In order to count integral points, we parameterize them by points on twisted universal torsors. These in turn are obtained through a Cox ring $R$: a universal torsor $Y$ can be described as a subvariety of $\Spec R$, see \cite[\S\,4]{FP16}.
The minimal desingularization $\St$ is obtained by a sequence of five blow-ups of $\PP^2$, and by~\cite{D14}, a Cox ring of $\St$ over $K$ is given by
\begin{equation*}
    R = K[a_1,\dots,a_9]/(a_1a_9+a_2a_8+a_3a_4^2a_5^3a_7);
\end{equation*}
it is a graded ring, and the degrees $(\deg(a_1),\dots,\deg(a_9))$ are
\begin{equation}\label{eq:degrees_a}
     (l_5,\ l_4,\ l_1-l_2,\ l_2-l_3,\ l_3,\ l_0-l_1-l_4-l_5,\ l_0-l_1-l_2-l_3,\ l_0-l_4,\ l_0-l_5)
\end{equation}
in $\Pic \St \cong \ZZ^6$ with standard basis $l_0,\dots,l_5$. Here, $l_0$ is the pull-back of the class of a line, and $l_1,\dots,l_5$ are the pull-backs of the classes of the exceptional curves of the five blow-ups.
This grading induces an action of the Néron--Severi torus $T = \widehat{\Pic}_\St$ on $\Spec R$, where an element $t\in T(K)$ acts on homogeneous elements $x\in R$ through multiplication with the character $\deg(x)$; this action restricts to $Y$. 

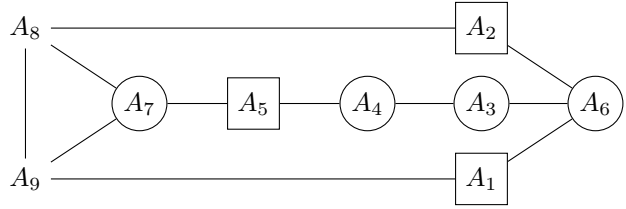
\begin{figure}[ht]
  \begin{center}
    \begin{tikzpicture}[square/.style={regular polygon,regular polygon sides=4}]
      \node (A9) at (0,-1) {$A_9$};
      \node (A8) at (0, 1) {$A_8$};
      \node[circle, draw, inner sep=.5ex] (A7) at (1.5, 0)  {$A_7$};
      \node[square, draw, inner sep=.2ex] (A5) at (3, 0)  {$A_5$};
      \node[circle, draw, inner sep=.5ex] (A4) at (4.5, 0)  {$A_4$};
      \node[circle, draw, inner sep=.5ex] (A3) at (6, 0)  {$A_3$};
      \node[circle, draw, inner sep=.5ex] (A6) at (7.5, 0) {$A_6$};
      \node[square, draw, inner sep=.2ex] (A2) at (6, 1)  {$A_2$};
      \node[square, draw, inner sep=.2ex] (A1) at (6,-1)  {$A_1$};

      \draw (A9) -- (A8);
      \draw (A9) -- (A7);
      \draw (A8) -- (A7);
      \draw (A8) -- (A2);
      \draw (A9) -- (A1);
      \draw (A6) -- (A2); 
      \draw (A6) -- (A1); 
      \draw (A7) -- (A5); 
      \draw (A5) -- (A4); 
      \draw (A4) -- (A3); 
      \draw (A3) -- (A6);
    \end{tikzpicture}
    \caption{The extended Dynkin diagram of the weak del Pezzo surface $\St$. The $(-1)$-curves are shown in squares and the $(-2)$-curves in circles.}\label{fig:dynkin-diagram}
  \end{center}
\end{figure}

Each generator $a_i$ corresponds to a divisor $A_i = V(a_i)$ defined in \emph{Cox coordinates}, with the first seven corresponding to negative curves. Concretely, $A_7$ is the exceptional divisor over $Q_1$, and $A_3+A_4+A_6$ is the exceptional divisor over $Q_2$, while $A_5$, $A_1$, and $A_2$ are the strict transforms of the lines $L$, $L'$, and $L''$, respectively.
The intersection properties of these divisors are encoded in the extended Dynkin diagram (Figure~\ref{fig:dynkin-diagram}), which also involves the two divisors $A_8$ and $A_9$ corresponding to the last two variables.

In particular, these divisors satisfy $[A_i] =\deg(a_i)\in \Pic \St$. We observe that $D=A_3+\dots+A_7$ has class $2l_0-l_1-\dots-l_5$, hence $-K_\St-D$ has class $l_0$.
For the following arguments, it will sometimes be convenient to work in the basis $([A_3],\dots,[A_8])$ of the Picard group.
In this basis, the degrees of the remaining generators are readily computed to be
\begin{align}
    [A_1]&=-[A_3]-[A_4]-[A_5]-[A_6]+[A_8],\label{eq:A1}\\ 
    [A_2]&=[A_3]+2[A_4]+3[A_5]+[A_7]-[A_8],\qquad \text{and}\label{eq:A2}\\
    [A_9]&=2[A_3]+3[A_4]+4[A_5]+[A_6]+[A_7]-[A_8].\nonumber
\end{align}

\subsection{Twisted universal torsors over $\OK$} Frei and Pieropan~\cite[\S\,4]{FP16} construct a model $\rho\colon \Yc\to\Sct$ of the universal torsor $Y \to \St$ as a subscheme of the spectrum of 
\begin{equation*}
    \Rc = \OK[a_1,\dots,a_9]/(a_1a_9+a_2a_8+a_3a_4^2a_5^3a_7).
\end{equation*}
Here, $\Sct$ is the model of $\St$ defined as a sequence of five blow-ups of $\PP^2$, whose centers are the closures of those of the sequence of five blow-ups defining $\St$.
This morphism is a torsor under the Néron--Severi torus $\Tc = \widehat{\Pic}_{\Sct/\OK}\cong \GGmOK^6$ dual to the relative Picard scheme. It is a universal torsor~\cite[Prop.~4.1]{FP16}, meaning that the class of every line bundle $\Lc$ in the relative Picard group,
corresponding to a character $\chi$ of $\Tc$, coincides with  $[\Yc\times^{\Tc} \AAA_{\OK}^1]\in \Pic_{\Sct/\OK}(\Sct) = \Pic \St$, where $\Tc$ acts on $\AAA^1$ through $\chi$ (note that this is the inverse of the usual type, owing to the Picard group's action on the spectrum of the Cox ring through the degree instead of inverse degree); for $d\in \Pic S$, write $\Lc_d$ for this push-forward bundle, and $L_d$ for the analogously defined push-forward of $Y$ over $K$. The global sections of these bundles are elements of $\Rc$ and $R$ of degree $d$, respectively.

Orbits of integral points on the twists of this torsor by elements of $H^1(\OK,\Tc) \cong \Cl_K^6$ parameterize integral points on $\Sct$, as for each integral point on $\Sct$, the fiber on precisely one of these twists is a trivial $\Tc$-torsor.
To describe $\Yc$ and its twists explicitly, we first need to introduce some more notation. For a $6$-tuple $\cfrb = (\cfr_3,\dots,\cfr_8)$ of nonzero fractional ideals of $\OK$, let $\Os_i = \cfrb^{[A_i]}$, where $[A_i]$ is interpreted as an element of $\ZZ^6$ via the basis $([A_3],\dots,[A_8])$, and exponentiation is to be read as the product of the component-wise powers. Let $\Os_{i*} = \Os_i^{\ne 0}$. We define the ideals $\afr_i = a_i\Os_i^{-1}$ for $a_i \in \Os_i$.

As found in the work of Frei and Pieropan~\cite[p.~777]{FP16}, the set of integral points on the \emph{twist $\crho\colon \cYc\to \Sct$ of $\Yc$ by $\cfrb$} is 
\begin{equation*}
    \cYc(\OK) = \bigwhere{(a_1,\dots,a_9) \in \Os_1\times \dots \times \Os_9}{&a_1a_9+a_2a_8+a_3a_4^2a_5^3a_7 = 0\\ &\text{\eqref{eq:coprimality} holds}},
\end{equation*}
involving the coprimality condition
\begin{equation}\label{eq:coprimality}
    \afr_i+\afr_j = \OK \text{ whenever $A_i$ and $A_j$ do not share an edge in Figure~\ref{fig:dynkin-diagram}},
\end{equation}
and an element
\begin{equation}\label{eq:torus_action}
    \text{$t\in \Tc(\OK) \cong (\OK^\times)^6$
acts on $\cYc(\OK)$ by multiplying each $a_i$ with $t^{[A_i]}$.}
\end{equation}
Its set of local integral points $\cYc(\Op)$ has an analogous shape, each $\Os_i$ being replaced by $\cfrb^{[A_i]} \Op$.

\subsection{Parametrizing integral points by orbits} In particular, letting $\Cs_0$ be a system of representatives of $\Cl_K$ consisting of integral ideals and containing $\OK$ as the principal ideal, we obtain a bijection
\begin{equation}\label{eq:bijection-projective}
    \coprod_{\cfrb \in \Cs_0^6} \cYc(\OK) / \Tc(\OK) \to \Sct(\OK).
\end{equation}

To restrict it to integral points away from $D$, write $\Ac_i = \overline{A_i}\subset \Sct$ and $\Dc = \Ac_3 + \dots + \Ac_7 = \overline {D}$ for the closures of the divisors $A_1,\dots,A_9$ and the boundary $D$, thus
obtaining a model $\Uct = \Sct \setminus \Dc$ of $\St \setminus D \cong U$. Defining $\cYcU = \cYc \times_{\Sct} \Uct$ and restricting~\eqref{eq:bijection-projective} yields a bijection
\begin{equation}\label{eq:bijection-open}
    \coprod_{\cfrb \in \Cs_0^6} \cYcU(\OK) / \Tc(\OK) \to \Uct(\OK).
\end{equation}

Let us describe the torsors on the left hand side more explicitly. The scheme $\cYcU$ is the complement of $V(a_3a_4a_5a_6a_7)$ in $\cYc$. Hence, its set of integral points is 
\begin{equation*}
    \cYcU(\OK) = \bigwhere{(a_1,\dots,a_9)\in \Os_1\times \dots \times \Os_9}
    {&a_1a_9+ a_2a_8+a_3a_4^2a_5^3a_7 = 0\\
    &\text{\eqref{eq:coprimality} holds},\\
    &\afr_3= \ldots =\afr_7=\OK};
\end{equation*}
the latter conditions are analogous to~\eqref{eq:coprimality} and can be obtained from~\cite[Theorem~2.7~(c)]{FP16}.
It implies $\cfr_3 = \dots = \cfr_7=\OK$ and $a_3,\dots, a_7\in \OK^\times$.
Moreover, it implies~\eqref{eq:coprimality}: the coprimality conditions involving at least one of $\afr_3,\dots,\afr_7$ follow immediately, and the remaining ones, only involving $\afr_1$, $\afr_2$, $\afr_8$, and $\afr_9$ follow from the torsor equation, which  now implies $\afr_1\afr_9+ \afr_2\afr_8 = \OK$. We thus obtain that $\cYcU(\OK)$ is
\begin{equation}\label{eq:integral_points_on_torsor}
    \{(a_1,\dots,a_9) \in \Os_1\times\Os_2\times \OK^\times \times \dots \times \OK^\times \times \Os_8\times \Os_9 : a_2a_8+a_1a_9+a_3a_4^2a_5^3a_7 = 0\}
\end{equation}
if $\cfr_3=\dots=\cfr_7=\OK$, and otherwise $\cYcU(\OK)=\emptyset$. This motivates the definition
\begin{equation}\label{eq:nonempty_classes}
    \Cs = \{(\OK,\dots,\OK,\cfr_8) : \cfr_8 \in \Cs_0\} \subset \Cs_0^6.
\end{equation}
Note that $\Os_1=\Os_8=\cfr_8$ and $\Os_2=\Os_9=\cfr_8^{-1}$ for $\cfrb \in \Cs$.

Having better understood the left hand side of~\eqref{eq:bijection-open}, we now turn to its right hand side and compare it to $\Uc(\OK)$.

\begin{lemma}\label{lem:local-points-coincide}
    The desingularization morphism induces an isomorphism $\Uct \to \Uc$. In particular, $\Uc(\OK)=\Uct(\OK)$.
\end{lemma}

\begin{proof}
    The composition $f$ of the torsor morphism $Y\to \St$ with the desingularization $\St \to S$ maps  $(a_1,\dots,a_9)\in \Yc(\Op)$ to
    \begin{equation}\label{eq:map-from-torsor}
       (a_2a_3a_4a_5a_6a_7a_8 : a_1^2a_2^2a_3^2a_4a_6^3: a_1a_2a_3^2a_4^2a_5^2a_6^2a_7 : a_1a_3a_4a_5a_6a_7a_9 : a_7a_8a_9).
    \end{equation}
    These five monomials, interpreted as global sections of the line bundle $\Lc_{-K_{\St}}$ on $\Sct$, do not vanish simultaneously: by~\cite[Prop.~4.1]{FP16}, each special fiber is a del Pezzo surface, and the above sections form a basis of the global anticanonical sections.
    They thus induce a morphism $\Sct\to \PP_{\OK}^4$ with image $\Sc$. Note that the monomials corresponding to $x_0$, $x_1$, and $x_3$ share the factors $a_3,\dots,a_7$, while the remaining factors $(a_2a_8,a_1a_2,a_1a_9)$ may not vanish simultaneously modulo any prime as a consequence of the comprimality conditions. It follows that the preimage of $\overline{L}=V_{\Sc}(x_0,x_1,x_3)$ is precisely $\Dc = V_{\Sct}(a_3\cdots a_7)$.

    We thus obtain a projective and birational restriction $\Uct\to \Uc$. This restriction is quasi-finite. Indeed, it only contracts the $(-2)$-curves above the singularity inside the generic fiber, and the situation is identical for the special fibers by~\cite[Prop.~4.1]{FP16}.

    Finally, computing derivatives of the defining equations, one may check that the singularities of all special fibers of $\Sc$ lie on the specializations of the two singularities outside $\Uc$ (including in characteristic $2$), so that the morphism $\Sc\to \Spec \OK$ is smooth, hence $\Sc$ is regular and in particular normal. We may thus conclude by Zariski's Main Theorem.
\end{proof}

\subsection{Lifting the height}\label{sec:height}

The elements $a_1a_9$, $a_2a_8$, and $a_1a_2a_3a_4a_5a_6$ of the Cox ring have log anticanonical degree, hence are global sections of the bundle $L_{[-K_{\St}-D]}$. As a consequence of the coprimality condition~\eqref{eq:coprimality}, they do not have a common zero and induce a metric on $L_{[-K_{\St}-D]}$, which in turn induces a log anticanonical height function with local factors
\begin{equation}\label{eq:local-factors-height}
    \Nt_v(x_1,\dots,x_9) = \max\{|x_2x_8|_v, |x_1x_2x_3x_4x_5x_6|_v, |x_1x_9|_v\}.
\end{equation}
on the universal torsor, which coincides with the height function $H$ described in the introduction~\eqref{eq:height-intro} in the following sense.

\begin{lemma}
    Let $\cfr\in \Cs$ and $(a_1,\dots,a_9)\in \cYcU(\OK)$. Then
    \begin{equation}\label{eq:height_on_torsor}
        H(\pi(\crho(a_1,\dots,a_9))) = \prod_{v\in \Omega_K}\Nt_v(a_1^{(v)},\dots,a_9^{(v)}).
    \end{equation}
\end{lemma}

\begin{proof}
    Using the description~\eqref{eq:map-from-torsor} that analogously describes all twists $\pi\circ \crho$, the left-hand side of~\eqref{eq:height_on_torsor} evaluates to 
    \begin{equation*}
        \prod_{v\in \Omega_K} \max\{|a_2a_3a_4a_5a_6a_7a_8|_v, |a_1a_2a_3^2a_4^2a_5^2a_6^2a_7|_v, |a_7a_8a_9|_v\},
    \end{equation*}
    which simplifies to the right-hand side after pulling out the common factor $|a_3\cdots a_7|_v$ from the maximum and using the product formula.
\end{proof}

The parameterization of integral points derived in this section can be summarized as in the following proposition.

\begin{prop}\label{prop:abstract_parameterization}
    The map $f\colon Y\to S$ given by \eqref{eq:map-from-torsor} induces a bijection between
    \begin{enumerate}
        \item the $\Tc(\OK)$-orbits under the action \eqref{eq:torus_action} on $\bigcup_{\cfrb \in \Cs} \cYcU(\OK)$ as in \eqref{eq:integral_points_on_torsor} and \eqref{eq:nonempty_classes}, and
        \item the set $\Uc(\OK)$ of integral points on our singular quartic del Pezzo surface $S$ defined by \eqref{eq:surface}.
    \end{enumerate}
    The height $H$ on $\Uc(\OK)$ as in \eqref{eq:height} lifts to $\cYcU(\OK)$ as in \eqref{eq:local-factors-height} and \eqref{eq:height_on_torsor}.
\end{prop}

\section{The expected asymptotic formula}\label{sec:expected}

\subsection{Tamagawa measures}\label{sec:tamagawa}

As a model of the universal torsor $Y\to \St$ under $T = \widehat{\Pic}_\St$, Frei and Pieropan construct a universal torsor $\Yc \to \Sct$ as a quasi-affine variety inside $\Spec \Rc$, where $\Rc\otimes_{\OK} K$ is the Cox ring $R$,
which the Néron--Severi torus $\Tc = \widehat{\Pic}_{\Sct/\OK}$ acts on.
If $L$ is a line bundle on $\St$, corresponding to a character $\chi\colon T \to \GGm$,
then $L$ corresponds to the $\GGm$-torsor $\chi_*Y$ in the sense that it is isomorphic to the contracted product $L \cong Y\times^T \AAA^1$, where $T$ acts on $\AAA^1$ through $\chi$, and analogously for the models. For an element $d$ of the Picard group, write $L_d$ and $\Lc_d$ for these two line bundles over $\St$ and $\Sct$, respectively.
Concretely, if $f\in R$ such that  $\St_f = \St \setminus V(f)$ (in Cox coordinates) is affine, then the set $H^0(\St_f,L_d)$ of sections is the set $R_f^{(d)}$ of elements of the localization of the Cox ring of degree $d$ (and making analogous definitions, $H^0(\Sct_f,\Lc_d) = \Rc_f^{(d)}$).

In order to compute the Tamagawa measures appearing in the expected formula, we are in need of a metric on the log anticanonical bundle $\omega_\St(D)^\vee$ (inducing one on its dual, the log canonical bundle), while the constructions in the previous sections furnish one on the bundle $L_{[-K_\St-D]}$ isomorphic to it. To obtain a norm on the former, we shall work with the isomorphism $\psi$ identifying the section $a_5a_7/a_1^2a_2^2 a_3 a_6^2$ of $L_{[K_\St+D]}$ with the section 
\begin{equation*}
    1_{D_3}1_{D_4}1_{D_5}1_{D_6}1_{D_7} \left(\ddd \frac{a_4a_5^2a_7}{a_1a_2a_6} \wedge \ddd\frac{a_8}{a_1a_3a_4a_5a_6} \right)
\end{equation*}
of $\omega_\St(D)$.
(As in~\cite[\S\,6]{DW24}, the section $1/a_1^2a_2^2 a_3^2 a_4a_6^3$ has degree $K_\St$;
it and the $2$-form in the latter expression are characterized, uniquely up to scalar, by the property of having neither zeroes nor poles on the open subvariety
$\St \setminus V(a_1a_2a_3a_4a_5a_6)$, coinciding with the affine open $D(f_1) = D(A_{7,8,9})$ in the notation of~\cite{FP16}).
\subsubsection*{Archimedean measures}

Over an archimedean place $v$, the volumes can be computed analogously to~\cite[Lem.~25]{DW24}, the only difference being that the parameters in the limits are possibly complex (which does not change the value of those limits) and that the renormalization factor $c_{\CC}^2=4\pi^2$ appears instead of $c_\RR^2=4$ when $v$ is complex~\cite[Eq.~(3.1), \S\,4.1]{CLT10}. It follows that 
\begin{equation}\label{eq:archimedean-comparison}
    \prod_{v \mid \infty} \tau_{{B_v},v}(Z_{B_v}(K_v)) = 4^{r_1}(4\pi^2)^{r_2} = \prod_{v\mid\infty} \omega_v,
\end{equation}
independently of $\Bb$.

\subsubsection*{Nonarchimedean measures}
 Let $\pfr$ be a nonarchimedean place and $\Sp = \Sct \times_{\OK} \Op$.
By flatness and as the special fiber is reduced and irreducible,
\begin{equation*}
    0=\Pic\Op \to \Pic\Sp \to \Pic\St\to 0
\end{equation*}
is exact, and it follows that $\Lc_{[K_\St+D]}|_{\Sp} \cong \omega_{\Sct/\OK}(\Dc)|_{\Sp}$ as their generic fibers coincide. The same sections as before form the unique primitive regular sections on $\Sct \setminus V(a_1a_2a_3a_4a_5a_6)\cong \AAA^2_{\OK}$ (note that the local chart spreads out) of the respective bundles that vanish nowhere on the generic fiber, up to units. It follows that there is an isomorphism between the respective bundles identifying these sections; that is, $\psi$ spreads out.

We note that the integrality condition together with $a_2a_8+a_1a_9\in \Op$ implies that~\eqref{eq:local-factors-height} is identically $1$ on $\Op$-integral points.
In particular, $\lVert1_D \omega\rVert_\pfr = 1$ for every primitive $2$-form; hence, $\lVert \cdot \rVert_\pfr$ is the model metric, and as the model $\Sct$ is smooth~\cite[Prop.~4.1]{FP16}, the equality
\begin{equation}\label{eq:tamagawa-vs-counting}
    \tau_{(\St,D),\pfr}(\Uct(\Os_\pfr)) =
    \frac{\# \Uct(\Os_\pfr/\pfr)}{\N \pfr^2},
\end{equation}
which a priori holds for almost all places, in fact holds for all of them. We are thus left to compute points on $\Uc$ modulo prime ideals.
To that end, note that $\tilde{\Sc}$ admits $\N \pfr^2 + 6 \N\pfr + 1$ points modulo $\pfr$, as a consequence of the Lefschetz trace formula (note that the higher cohomology groups of $\Os_{\St}$ are trivial by Kawamata--Viehweg vanishing), the Picard group of $\Sct$ being split of rank $6$, and Poincaré duality.
Each of the exceptional curves $A_3,\dots,A_7$ is isomorphic to a copy of $\PP^1$,
so after taking their four intersection points into account, $\Dc$ admits $5\N \pfr+1$ points modulo $\pfr$.
Thus $\Uc(\Os_\pfr/\pfr) = \N\pfr^2 + \N\pfr$, and Equation~\eqref{eq:tamagawa-vs-counting} together with Lemma~\ref{lem:local-points-coincide} yields
\begin{equation*}
    \tau_{(\St,D),\pfr}(\Uc(\Os_\pfr)) = 1+ \frac{1}{\N\pfr}.
\end{equation*}
The Picard group of $U$ is split of rank $1$, hence we multiply these volumes with the local factors $\lambda_\pfr = 1-\N\pfr^{-1}$ of its associated Artin $L$-function $\zeta_K$ as convergence factors, yielding the local densities
\begin{equation*}
    \lambda_\pfr\tau_{(\St,D),\pfr}(\Uc(\Os_\pfr)) = 1- \frac{1}{\N\pfr^2} = \omega_\pfr
\end{equation*}
as in~\eqref{eq:omega_v}
and the factor $\rho_K$ in~\eqref{eq:c} as the residue of~$\zeta_K$.

In total, 
\begin{equation}\label{eq:c_fin-comparison}
    c_\fin = \rho_K \prod_{\pfr} \omega_{\pfr}.
\end{equation}

\subsection{The effective cone constant $\alpha$: Solving the jigsaw puzzle}\label{sec:alpha}

The analytic Clemens complex has $4^{r_1+r_2}=4^{q+1}$ maximal faces, each with an associated Picard group, effective cone and polytope whose volume appears in the expected formula. As the actual asymptotic formula only involves one such volume, the key observation at this point will be that these $4^{q+1}$ polytopes fit together to one larger polytope --- as they are defined in the duals of different Picard groups, this involves working in specific, compatible presentations of these groups.

\subsubsection*{Picard groups and gluing}

We begin the computation of the effective constants with general observations on their definition. For a face $B$ of the Clemens complex $\Cs_v(D)$ at a place $v$, write
\begin{equation*}
    U_B=\St\setminus \bigcup_{\substack{i\in \{3,\dots,7\},\\ A_i\nprec B}} A_i
\end{equation*}
for the complement of the components of $D$ not belonging to $B$.
The variants of the Picard group pertinent to the study of integral points can be described as the fiber products
\begin{equation*}
    \Pic(U;\Bb) = \sideset{}{_{\Pic U}}\prod_{v\mid \infty} \Pic U_{B_v}
\end{equation*}
associated with each face $\Bb = (B_v)_{v\mid\infty}\in \Cs(D)$~\cite[p.~577]{Wil24}.
For any open subvariety $W$ of $\St$, write $\sim_W$ for rational equivalence on $W$ and $[\cdot]_W$ for the resulting equivalence classes whenever this is necessary to avoid confusion between the various relations.
The vector spaces $\Pic(U;\Bb)_\RR$ contain effective cones $\Eff(U;\Bb)$.
We will intersect their duals $\Eff(U;\Bb)^\vee$ each with a certain half-space; the volumes of the resulting polytopes $P_{\Bb}$ for maximal faces are part of the leading constant~\cite[Rmk.~2.2.9~(iv)]{Wil24}.
If $\Bb'\prec \Bb$ is a subface, the projection (or pullback) map $\Pic(U;\Bb)\to \Pic(U;\Bb')$ induced by the pullbacks $\Pic U_{B_v} \to \Pic U_{B'_v}$ corresponds to an injection of the duals, including an injection $\Eff(U;\Bb')^\vee\to \Eff(U;\Bb)^\vee$.
Collectively, these injections provide gluing data for the polytopes $P_{\Bb}$.

Working in ``compatible'' bases, these injections simply become inclusions of subsets. More formally, these bases for varying $\Bb$ will be the images of a single subset $\mathcal{B}$ of the fiber product of $\#\Omega_\infty$ copies of $\Pic X$ over $\Pic U$ in its various quotients $\Pic(U;\Bb)$; they thus induce isomorphisms between the duals of the latter to $\langle \mathcal{B}\rangle^\vee$. It will moreover turn out that, at least for this particular variety, the interiors of the polytopes are disjoint, so that gluing them can be achieved by simply taking their union.
It would be interesting to know whether this can be achieved in general.

\subsubsection*{Bases}

To begin with, observe that $\Pic U = \ZZ [A_8]_U$ is of rank one, while $\Pic U_B = \ZZ[A_1]_{U_B} \oplus \ZZ[A_2]_{U_B}\oplus \ZZ[A_8]_{U_B}$ for every $B\in \Cmax_v(D)$.
In particular, the fiber product $\Pic(U;\Bb)$ is of rank $2(q+1)+1$.
In order to obtain a basis, order the archimedean places $\Omega_\infty = \{v_0,\dots,v_q\}$, and fix a maximal face $\Bb = (B_v)_{v\mid\infty}$ for now. Write $B_n = B_{v_n}$ and $U_n = U_{B_n}$ to simplify notation.
As $A_1 \sim_U A_8$ and $A_2 \sim_U -A_8$ by \eqref{eq:degrees_a}, working modulo $[A_3]_\St, \dots, [A_7]_\St$, the classes $e_{n,1} = [A_1-A_8]_{U_n}$ and $e_{n,2} = [A_2 + A_8]_{U_n}$ become trivial on $U$. 
It follows that the group $\Pic(U;\Bb)$ is freely generated by the divisor class
$E_0=([A_8])_{v\mid\infty}$ and the $2(q+1)$ classes $E_{n,i} = (0,\dots,0,e_{n,i},0,\dots, 0)$, where $i\in \{1,2\}$, $n\in \{0,\dots,q\}$, and the nontrivial entry is the $(n+1)$th one, corresponding to the place $v_n$.

Write
\begin{equation*}
    \ab = (a_0,(a_{0,1}, a_{1,1}),\dots, (a_{q,1},a_{q,2})) = a_0E_0 + \sum a_{n,i} E_{n,i}
\end{equation*}
in this basis.

\subsubsection*{Cones and polytopes}
Our next task is to describe the effective cone $\Eff(U;\Bb)$ in $\Pic(U;\Bb)_\RR$ in this basis. As the effective cone of $\St$ is generated by the classes of $A_1,\dots,A_7$, that of $\Pic(U;\Bb)$ is generated by 
\begin{align*}
    ([A_1]_{U_0},\dots,[A_1]_{U_q}) &= (1, (1,0),\dots, (1,0)),\\
    ([A_2]_{U_0},\dots,[A_2]_{U_q}) &= (-1, (0,1),\dots, (0,1)),
\end{align*}
and the $2q+2$ classes 
\begin{equation*}
    (0,\dots,0,[A_i]_{U_n},0,\dots,0),
\end{equation*}
where $0\le n\le q$ and $A_i\prec B_n$ (as $[A_i]_{U_n}=0$ whenever $A_i\nprec B_n$). 
The identities~\eqref{eq:A1} and~\eqref{eq:A2} yield the relations 
\begin{equation*}
    A_3 + A_4 + A_5 + A_6 \sim_{\St}  A_8 -A_1 \quad\text{and}\quad A_3 + 2A_4 + 3A_5 + A_7 \sim_{\St} A_2+A_8,
\end{equation*}
and 
together with $A_i\sim_{U_n} 0$ if $A_i\nprec B_n$,
we readily obtain
\begin{align*}
    [A_5]_{U_n} &= - e_{n,1},        & [A_7]_{U_n}&=3e_{n,1} + e_{n,2} & \text{if } B_n &= (57),            \\
    [A_4]_{U_n} &= -3e_{n,1}-e_{n,2},& [A_5]_{U_n}&=2e_{n,1} + e_{n,2} & \text{if } B_n &= (45),            \\
    [A_3]_{U_n} &= -2e_{n,1}-e_{n,2},& [A_4]_{U_n}&= e_{n,1} + e_{n,2} & \text{if } B_n &= (34), \text{ and}\\
    [A_3]_{U_n} &=  e_{n,2},         & [A_6]_{U_n}&=-e_{n,1} - e_{n,2} & \text{if } B_n &= (36).              
\end{align*}
Moreover, as $-K = A_2+\dots+A_8$, the log anticanonical class has the representation
\begin{equation*}
    - \Kb_D = ([-K-D]_{U_0},\dots,[-K-D]_{U_q}) = (0, (0,1),\dots,(0,1)).
\end{equation*}
Working with vectors $\ab = (a_0,(a_{0,1}, a_{1,1}),\dots (a_{q,1},a_{q,2}))$ in $\Pic(U;\Bb)_\RR^\vee$ expressed in the dual basis,
the polytope
\begin{equation*}
    P_{\Bb} = \{ E \in \Eff(U;\Bb)^\vee : \langle -\Kb_D, E\rangle\le 1 \}
\end{equation*}
associated with the face $\Bb$~\cite[Rmk.~2.2.9~(iv)]{Wil24} is thus described by the three inequalities
\begin{equation*}
    a_0 + \sum_{n=0}^q a_{n,1} \ge 0,\quad -a_0 + \sum_{n=0}^q a_{n,2}\ge 0,\quad\text{and}\quad \sum_{n=0}^q a_{n,2}\le 1
\end{equation*}
that are independent of the concrete face $\Bb$ and for each $0\le n \le q$, two more inequalities in only the two variables $a_{n,1}$ and $a_{n,2}$, which depend only on $B_n$.
As the face $B_n$ varies over its four possible values, these inequalities are
\begin{equation}\label{eq:cone-ineqs}
   \begin{aligned}
        - a_{n,1}        &\ge 0,& 3a_{n,1} + a_{n,2}&\ge 0 &\text{if } B_n &= (57),            \\    
        -3a_{n,1}-a_{n,2}&\ge 0,& 2a_{n,1} + a_{n,2}&\ge 0 &\text{if } B_n &= (45),            \\
        -2a_{n,1}-a_{n,2}&\ge 0,&  a_{n,1} + a_{n,2}&\ge 0 &\text{if } B_n &= (34), \text{ and}\\
          a_{n,2}        &\ge 0,& -a_{n,1} - a_{n,2}&\ge 0 &\text{if } B_n &= (36).                       
\end{aligned} 
\end{equation}

\subsubsection*{Taking unions}

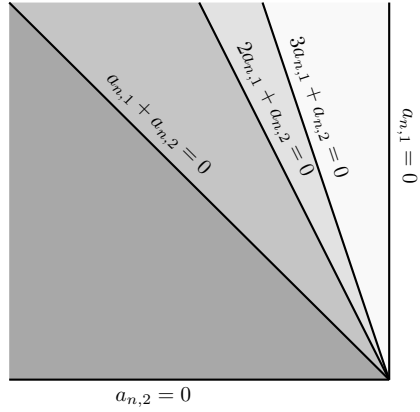
\begin{figure}[ht]
\begin{tikzpicture}
  \def\pcolora{gray!90}
  \def\pcolorb{gray!20}
  \def\popacity{0.4}
  \path[fill=\pcolora, opacity=0.05] (0,0) -- (0,   5) -- (-5/3,5);
  \path[fill=\pcolora, opacity=0.25] (0,0) -- (-5/3,5) -- (-5/2,5);
  \path[fill=\pcolora, opacity=0.5] (0,0) -- (-5/2,5) -- (-5,  5);
  \path[fill=\pcolora, opacity=0.75] (0,0) -- (-5,  5) -- (-5,  0);
  \begin{scope}[every node/.style={scale=.8}]
  \draw[thick] (0,0) --node[above left, sloped,rotate=180]{$a_{n,1}=0$} (0,   5);
  \draw[thick] (0,0) --node[above left, sloped]{$3a_{n,1}+a_{n,2}=0$}   (-5/3,5);
  \draw[thick] (0,0) --node[above left, sloped]{$2a_{n,1}+a_{n,2}=0$}   (-5/2,5);
  \draw[thick] (0,0) --node[above left, sloped]{$a_{n,1} + a_{n,2}=0$}  (-5,  5);
  \draw[thick] (0,0) --node[below left, sloped]{$a_{n,2}=0$}            (-5,  0);
\end{scope}
\end{tikzpicture}
\caption{The cones in the variables $a_{n,1}$ and $a_{n,2}$ described by the four sets of inequalities~\eqref{eq:cone-ineqs}, ordered counterclockwise. The actual polytopes are finite and described by additionally delimiting lines that depend on the other variables but not on the face $B_n$.}\label{fig:cones}
\end{figure}

Observe that when passing from one line (I) to the next or the previous (II), one inequality gets reversed (hence, the interiors of these polytopes are disjoint) and the new inequality in (II) is a positive linear combination of the ones in (I); in particular, this new inequality is implied by the ones in (I). We may thus take the union of the polytopes described by (I) and (II) by deleting the inequality that gets reversed and keeping the two others. Repeating this process for all three such steps leaves us with the two extremal inequalities $a_{n,1}\le 0$ and $a_{n,2}\ge 0$ (Figure~\ref{fig:cones}), while leaving all remaining inequalities, which are independent of $B_n$, intact.

Finally, repeating this process for all $n\in \{0,\dots,q\}$, we arrive at the disjoint union
\begin{equation*}
    P = \coprod_{\Bb\in \Cmax(D)} P_{\Bb} =
    \bigwhere{
        \ab \in \RR^{2q+3}
        }{
        \substack{
            a_0 + \sum_{n=0}^q a_{n,1} \ge 0,\quad  -a_0 + \sum_{n=0}^q a_{n,2}\ge 0,\\
            \sum_{n=0}^q a_{n,2}\le 1,\\
            a_{n,1}\le 0,\ a_{n,2}\ge 0 \text{ for all } n\in\{0,\dots,q\}
        }
        }.
\end{equation*}

\subsubsection*{Comparison with the main theorem}
To compare its volume with that appearing in Theorem~\ref{thm:main_concrete}, a change of variables flipping the signs of the coordinates $a_{n,1}$ and then replacing the two coordinates $a_{0,i}$ by $a_i=\sum_{n=0}^q a_{n,i}$  is convenient; note that both coordinates are nonnegative. As this change of variables has determinant $\pm 1$, the polytope $P$ is isometric to
\begin{equation*}
    P' = \bigwhere{
        (a_0,a_1,a_2,(a_{1,1},a_{1,2}),\dots,(a_{q,1},a_{q,2})) \in \RR^{2q+3}
        }{
        \substack{
            0\le a_1\le a_0 \le  a_2,\
            0\le a_2\le 1,\\
            \sum_{n=1}^q a_{n,1} \le a_1,\\ \sum_{n=1}^q a_{n,2} \le a_2, \\
            a_{n,1}\ge 0,\ a_{n,2}\ge 0\\ \text{ for all } n\in\{1,\dots,q\}
        }
        }.
\end{equation*}
This polytope is a pyramid with the origin as its apex and whose base $P'_0$ is the intersection of $P'$ with $a_2=1$, that is,
\begin{equation*}
    P'_0 = \bigwhere{
        (a_0,a_1,(a_{1,1},a_{1,2}),\dots,(a_{q,1},a_{q,2})) \in \RR^{2q+2}
        }{
        \substack{
            0\le a_1 \le a_0 \le 1,\\
            \sum_{n=1}^q a_{n,1} \le a_1,\\ \sum_{n=1}^q a_{n,2} \le 1, \\
            a_{n,1}\ge 0,\ a_{n,2}\ge 0\\ \text{ for all } n\in\{1,\dots,q\}
        }
        },
\end{equation*}
and $\vol(P') = \vol(P'_0)/(2q+3)$. For $a \in \RR_{\ge 0}$, note that
\begin{equation*}
    \vol\left\{(s_1,\dots,s_q)\in \RR^q_{\ge 0} : \sum_{n=1}^q s_n\le a\right\} = \frac{a^q}{q!},
\end{equation*}
so that
\begin{equation*}
    \vol{P_0'} = \frac{1}{q!^2} \int_{0\le a_1\le a_0\le 1} a_1^q \ddd a_0 \ddd a_1,
\end{equation*}
which coincides with $\alpha$ by its definition~\eqref{eq:alpha} upon the change of variables $(a_0,a_1) = (1-t_1,1-t_1-t_2)$ of determinant $1$.

This proves Proposition~\ref{prop:puzzle}. Moreover, together with~\eqref{eq:archimedean-comparison}, this implies that
\begin{equation}\label{eq:c_infty-comparison}
    c_\infty = \frac{\alpha}{\Delta_K} \prod_{v\mid\infty} \omega_v.
\end{equation}
\begin{proof}[Deduction of Theorem~\ref{thm:counting-expected-formula} from Theorem~\ref{thm:main_concrete}]
    Equations~\eqref{eq:c},~\eqref{eq:c_fin-comparison}, and~\eqref{eq:c_infty-comparison} imply that $c = c_\fin c_\infty$, while $\rk \Pic U=1$ and $\dim_v(\Cs_v(D))=1$ for all $v$ (the maximal faces being $1$-dimensional edges), so that $b=1 +2(q+1)$, there being $q+1$ archimedean places.
\end{proof}

\subsubsection*{Example: $q=1$} 
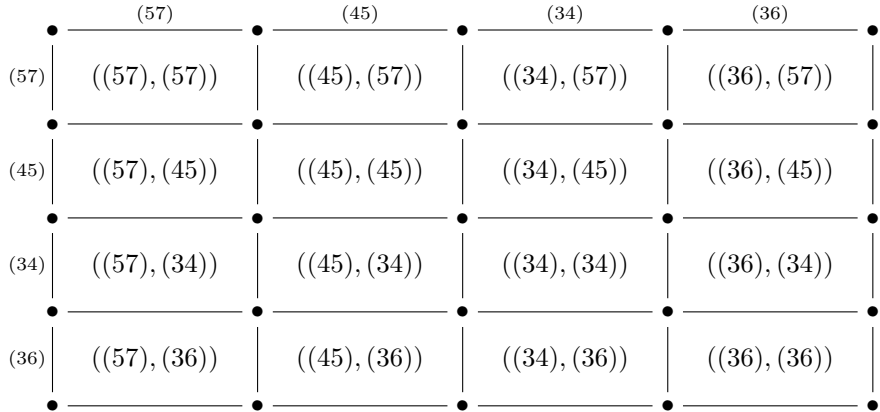
\begin{figure}[ht]
\[\begin{tikzcd}[sep=tiny,cramped]
	\bullet && \bullet && \bullet && \bullet && \bullet \\
	& {((57),(57))} && {((45),(57))} && {((34),(57))} && {((36),(57))} \\
	\bullet && \bullet && \bullet && \bullet && \bullet \\
	& {((57),(45))} && {((45),(45))} && {((34),(45))} && {((36),(45))} \\
	\bullet && \bullet && \bullet && \bullet && \bullet \\
	& {((57),(34))} && {((45),(34))} && {((34),(34))} && {((36),(34))} \\
	\bullet && \bullet && \bullet && \bullet && \bullet \\
	& {((57),(36))} && {((45),(36))} && {((34),(36))} && {((36),(36))} \\
	\bullet && \bullet && \bullet && \bullet && \bullet
	\arrow["{(57)}", no head, from=1-1, to=1-3]
	\arrow["{(57)}"', no head, from=1-1, to=3-1]
	\arrow["{(45)}", no head, from=1-3, to=1-5]
	\arrow[no head, from=1-3, to=3-3]
	\arrow["{(34)}", no head, from=1-5, to=1-7]
	\arrow[no head, from=1-5, to=3-5]
	\arrow["{(36)}", no head, from=1-7, to=1-9]
	\arrow[no head, from=1-7, to=3-7]
	\arrow[no head, from=1-9, to=3-9]
	\arrow[no head, from=3-1, to=3-3]
	\arrow[no head, from=3-3, to=3-5]
	\arrow[no head, from=3-3, to=5-3]
	\arrow[no head, from=3-5, to=3-7]
	\arrow[no head, from=3-5, to=5-5]
	\arrow[no head, from=3-7, to=3-9]
	\arrow[no head, from=3-7, to=5-7]
	\arrow[no head, from=3-9, to=5-9]
	\arrow["{(45)}", no head, from=5-1, to=3-1]
	\arrow[no head, from=5-1, to=5-3]
	\arrow[no head, from=5-3, to=5-5]
	\arrow[no head, from=5-3, to=7-3]
	\arrow[no head, from=5-5, to=5-7]
	\arrow[no head, from=5-5, to=7-5]
	\arrow[no head, from=5-7, to=5-9]
	\arrow[no head, from=5-7, to=7-7]
	\arrow[no head, from=5-9, to=7-9]
	\arrow["{(34)}", no head, from=7-1, to=5-1]
	\arrow["{(36)}"', no head, from=7-1, to=9-1]
	\arrow[no head, from=7-3, to=7-1]
	\arrow[no head, from=7-3, to=7-5]
	\arrow[no head, from=7-3, to=9-3]
	\arrow[no head, from=7-5, to=7-7]
	\arrow[no head, from=7-5, to=9-5]
	\arrow[no head, from=7-7, to=7-9]
	\arrow[no head, from=7-7, to=9-7]
	\arrow[no head, from=7-9, to=9-9]
	\arrow[no head, from=9-1, to=9-3]
	\arrow[no head, from=9-3, to=9-5]
	\arrow[no head, from=9-5, to=9-7]
	\arrow[no head, from=9-7, to=9-9]
\end{tikzcd}\]
\caption{The product $\Cs(D) = \prod_v \Cs_v(D)$ if $q=1$.}
\label{fig:Clemens-product}
\end{figure}

When $K$ has precisely two archimedean places, there are $16$ maximal faces of the Clemens complex (Figure~\ref{fig:Clemens-product}), each with a corresponding polytope.
The polytope associated with the face $((57),(57))$ has empty interior as a consequence of the analytic obstruction: the subvariety $U_{(57),(57)} = \St \setminus (A_3\cup A_3\cup A_6)$ admits a nonconstant regular function that makes the corresponding effective cone contain a line, so that its dual is not of full dimension.

The base $P_0'$ of the pyramid that is their union is defined inside a four-dimensional space with coordinates $(a_0,a_1,a_{1,1},a_{1,2})$. The intersection of this polytope with different planes by fixing $a_0$ and $a_1$ only depends on the latter coordinate, and the picture changes qualitatively when $a_1$ crosses the values $1/3$ and $1/2$: some polytopes whose intersection with the plane have full dimension in one range cease to do so at these values.

Figure~\ref{fig:jigsaw} thus depicts three such cross-sections of the polytope $P_0'$: with a plane with $a_1=0.2$ (left), $a_1=0.4$ (middle), or $a_1 =0.6$ (right), and arbitrary fixed $a_0$. The hue represents $B_0$ and opacity represents $B_1$. For instance, the bright orange polygon to the top left is associated with $((36),(57))$, and the dark green one on the bottom right is associated with $((57),(36))$.
Each of the 7, 11, and 11 (respectively) smaller polytopes is a triangle or quadrilateral, and in each case, they fit together to form an $(a_1\times 1)$-rectangle.

\section{Counting I: Preliminary reductions}\label{sec:parameterization}

We now begin the proof of Theorem~\ref{thm:main_concrete}. The first step is a reformulation of the parameterization of integral points obtained in Section~\ref{sec:torsor}.

\subsection{Parameterization of integral points}

We identify $\Tc(\OK)$ with $(\OK^\times)^6$ via the standard basis $\ell_0,\dots,\ell_5$ of $\Pic\St$. Let $\Fs$ be a fundamental domain for the action of the subgroup $U_K \times (\OK^\times)^5 \subset (\OK^\times)^6$ on $Y(K) \cap ((K^\times)^7 \times K^2)$. We will construct $\Fs$ explicitly in Section~\ref{sec:fund_domain}. We obtain the following reformulation of Proposition~\ref{prop:abstract_parameterization}.

\begin{prop}\label{prop:parameterization}
    For $\cfr \in \Cs$, let $M_\cfr(B)$ be the set of all
    \begin{equation*}
        (a_1,a_2,a_3,a_4,a_5,a_6,a_7,a_8,a_9) \in \Os_1\times\dots\times\Os_{9}
    \end{equation*}
    satisfying the torsor equation
    \begin{equation}\label{eq:torsor}
            a_2a_8+a_1a_9+a_3a_4^2a_5^3a_7 = 0,
    \end{equation}
    the height condition
    \begin{equation}\label{eq:height}
        \prod_{v \mid \infty} \Nt_v(a_1^{(v)},\dots,a_9^{(v)})\le B,
    \end{equation}
    the fundamental domain condition
    \begin{equation}\label{eq:fundamental_domain}
        (a_1,\dots,a_9) \in \Fs \subset (K^\times)^7 \times K^2,
    \end{equation}
    and the integrality condition
    \begin{equation*}
        a_3,a_4,a_5,a_6,a_7 \in \OK^\times.
    \end{equation*}

    Then
    \begin{equation*}
        N_{\Uc,V,H}(B) = \frac{1}{|\mu_K|} \sum_{\cfrb \in \Cs}|M_\cfrb(B)|.
    \end{equation*}
\end{prop}

\begin{proof}
    This follows from Proposition~\ref{prop:abstract_parameterization}. Since $\Vt = \pi^{-1}(V)$ is the complement of the negative curves $A_1,\dots,A_7$, the integral points in the complement $V$ of the three lines on $S$ lift to points with $a_1\cdots a_7\ne 0$ on the torsor. For the height function, we observe that
    \begin{equation*}
        \prod_{\pfr} \Nt_\pfr(a_1,\dots,a_9) = \N(\cfrb^{[-K_\St-D]})
    \end{equation*}
    for all $(a_1,\dots,a_9)\in \cYc(\OK)$. Note that $\N(\cfrb^{[-K_\St-D]})=1$ whenever $\cYcU(\OK)\ne \emptyset$. The factor $|\mu_K|^{-1}$ appears since $U_K$ has index $|\mu_K|$ in $\OK^\times$.
\end{proof}

\subsection{Construction of a fundamental domain} \label{sec:fund_domain}

In our fundamental domain, we will make sure that $a_1,\dots,a_5$ as well as the factors in the height function have `balanced' embeddings. In particular, this will force $a_3=a_4=a_5=1$.

We define $F(\infty), F(B), \Fs_1$ as in \cite[\S\,5]{FP16}. Let $\ab'=(a_1,\dots,a_5) \in (K^\times)^5$. For $(x_{6v},x_{7v},x_{8v}) \in K_v^3$, let
\begin{equation*}
    \Nt_v(\ab';x_{6v},x_{7v},x_{8v}) = \Nt_v\left(a_1^{(v)},\dots,a_5^{(v)},x_{6v},x_{7v},x_{8v},-\tfrac{a_2^{(v)}x_{8v}+a_3^{(v)}a_{4}^{(v)2}a_5^{(v)3}x_{7v}}{a_{1}^{(v)}}\right).
\end{equation*}
For $(x_{6v},x_{7v})_v \in \prod_{v \mid \infty} (K_v^\times)^2$, let
\begin{equation*}
    S_F(\ab',(x_{6v},x_{7v})_v;\infty) = \left\{(x_{8v})_v \in \prod_{v \mid \infty} K_v: (\log \Nt_v(\ab';x_{6v},x_{7v},x_{8v}))_v \in F(\infty)\right\}.
\end{equation*}
For $a_6,a_7 \in K^\times$, let
\begin{equation*}
    \Fs_0(\ab',a_6,a_7) \coloneq \{a_8 \in K : \sigma(a_8) \in S_F(\ab',(a_6^{(v)},a_7^{(v)})_v;\infty)\}.
\end{equation*}
Then the set of all $(a_1,\dots,a_9) \in (K^\times)^7 \times K^2$ with $\ab' \in \Fs_1^5$, $a_8 \in \Fs_0(\ab',a_6,a_7)$ and $a_9$ satisfying \eqref{eq:torsor} is a fundamental domain $\Fs$ as in Proposition~\ref{prop:parameterization}.

Let $S_F(\ab',(x_{6v},x_{7v})_v;B)$ be defined as
\begin{equation*}
     \left\{(x_{8v})_v \in \prod_{v \mid \infty} K_v : (\log \Nt_v(\ab';x_{6v},x_{7v},x_{8v}))_v \in F(B^{1/d})\right\}
\end{equation*}
and
\begin{equation*}
    \Fs_0(\ab',a_6,a_7;B) \coloneq \{a_8 \in K : \sigma(a_8) \in S_F(\ab',(a_6^{(v)},a_7^{(v)})_v;B)\}.
\end{equation*}
Then $(a_1,\dots,a_9) \in (K^\times)^7 \times K^2$ with \eqref{eq:torsor} satisfies \eqref{eq:height} and \eqref{eq:fundamental_domain} if and only if $\ab' \in \Fs_1^5$ and $a_8 \in \Fs_0(\ab',a_6,a_7;B)$.

\subsection{The dependent variable}\label{sec:dependent_variable}

As a first step towards our counting problem, we eliminate the variable $a_9$, introducing a congruence condition on $a_8$, thus rephrasing the problem as a lattice point count.

Let $\ab' \in \Os_*'$ and $a_6,a_7 \in \OK^\times$. If a solution to \eqref{eq:torsor} exists, we must have $\afr_1+\afr_2=\OK$; let $\theta_0(\afrb')$ be the indicator function of this coprimality condition. Conversely, suppose that $\theta_0(\afrb') = 1$. If

\begin{equation}\label{eq:congruence}
    \congr{a_2a_8}{-a_3a_4^2a_5^3a_7}{a_1\Os_9},
\end{equation}
then we obtain a unique
\begin{equation*}
    a_9=-(a_2a_8+a_3a_4^2a_5^3a_7)/a_1 \in \Os_9
\end{equation*}
satisfying \eqref{eq:torsor}; otherwise, no such $a_9$ exists.

The congruence condition \eqref{eq:congruence} can be rewritten as $a_8 \in \Gs = \Gs(\ab',a_6,a_7,\cfrb)$ for a shifted lattice $\Gs = \gamma + \Gs'$ with $\Gs' = \afr_1\Os_8$ (which we regard as an additive subgroup of $K$) and some $\gamma = \gamma(\ab',a_6,a_7,\cfrb) \in K$. Indeed, since $a_1\Os_9+a_2\Os_8=\afr_1+\afr_2=\OK$, there is always a solution $a_8$ of \eqref{eq:congruence}, and any two solutions differ by a solution $d \in \Os_8$ to $a_2d \in \afr_1$ which is equivalent to $\afr_1 \mid (a_2d)=\afr_2d\Os_2$ and hence to $\afr_1 \mid d\Os_2$ and hence to $\afr_1\Os_8 \mid d$, i.e., $d \in \afr_1\Os_8$. Hence
\begin{equation}\label{eq:lattice_points}
    |\{(a_8,a_9) \in \Os_8 \times \Os_9 : \eqref{eq:torsor}, \eqref{eq:height}, \eqref{eq:fundamental_domain}\}| = |\Gs(\ab',a_6,a_7,\cfrb) \cap \Fs_0(\ab',a_6,a_7;B)|.
\end{equation}

\subsection{Restricting $a_6,a_7$}\label{sec:restrictions}

In order to smoothly count lattice points in $\Gs'$ in the next subsection, we introduce some preliminary truncations that allow to control the error terms arising in that count. 

By the construction of our fundamental domain, the height condition \eqref{eq:height} implies that
\begin{equation}\label{eq:balanced_heights}
\Nt_v(a_1^{(v)},\dots,a_9^{(v)}) \ll B
\end{equation}
for all $v \mid \infty$.

Motivated by $\deg(a_6)=\ell_0-\ell_1-\ell_4-\ell_5$ and $\deg(a_1\cdots a_5)=\ell_1+\ell_4+\ell_5$ (and similar observations for $a_7, a_8$), we define
\begin{equation}\label{eq:B6B7B8}
    B_{6v}=\frac{B^{1/d}}{a_1^{(v)}\cdots a_5^{(v)}},\quad B_{7v}=\frac{B^{1/d}}{a_3^{(v)}a_4^{(v)2}a_5^{(v)3}},\quad B_{8v} = \frac{B^{1/d}}{a_2^{(v)}}
\end{equation}
and write $B_i \coloneq \prod_{v \mid \infty} |B_{iv}|_v$ for $i=6,7,8$, so that we have
\begin{equation}\label{eq:bound_ai}
    |a_i|_v \ll B_{iv} \asymp B_i^{d_v/d}.
\end{equation}
Moreover, we certainly have $\N(\afr_1\afr_2)\le B$.

We want to bootstrap these estimates to the conditions
\begin{equation}\label{eq:restrict_a6_a7}
    |a_i|_v \le (B_i/T)^{d_v/d}
\end{equation}
for $i=6,7$, as well as
\begin{equation}\label{eq:restrict_a1a2}
    \N(\afr_1\afr_2) \le \frac{B}{T}
\end{equation}
for a parameter $T=(\log B)^d$. Let
\begin{equation*}
    \Os'_* = \Os_{1*}\times \Os_{2*} \times \OK^\times \times \OK^\times \times \OK^\times.
\end{equation*}

\begin{prop}\label{prop:restrict_a6_a7}
    We have
    \begin{align*}
        |M_\cfrb(B)| ={}&\sums{\ab' \in \Os'_* \cap \Fs_1^5\\ \eqref{eq:restrict_a1a2}} \sums{a_6,a_7 \in \OK^\times\\\eqref{eq:restrict_a6_a7}} |\{(a_8,a_9) \in \Os_8 \times \Os_9 : \eqref{eq:torsor}, \eqref{eq:height}, \eqref{eq:fundamental_domain}\}|\\
        &+ O(B(\log B)^{1+2q}\log\log B).
    \end{align*}
\end{prop}

\begin{proof}
    We need to discard tuples with $|a_i|_v>(B_i/T)^{d_v/d}$ for $i=6$ or $i=7$ and some $v$ or with $B/T<\N(\afr_1\afr_2) \le B$. In view of \eqref{eq:balanced_heights}, this restricts one of $\max_v |a_6|_v$,  $\max_v |a_7|_v$ and $\N(\afr_1\afr_2)$ to one of $O(\log T)=O(\log\log B)$ many dyadic intervals.

    Since $a_8$ satisfies a congruence condition modulo $\Gs'=\afr_1\Os_8$ and $a_{8v}$ is bounded by $B_{8v}$, the number of solutions $(a_8,a_9)$ for a given tuple of the remaining variables can be estimated using \cite[Lemma~7.1]{FP16} as
    \begin{equation*}
        \ll \frac{B_8}{\N(\afr_1)\Os_8)} +1 \ll \frac{B}{\N(\afr_1\afr_2)}.
    \end{equation*}
    We can thus bound the number of such tuples by
    \[\sum_{a_1,a_2} \sum_{a_6,a_7 \in \Os_K} \frac{B}{\N(\afr_1\afr_2)}.\]
    Here, by \cite[Lemma~3.1]{BD25integral}, the number of $a_6$ and $a_7$ is $O((\log B)^{q})$ each. In view of \cite[Lemma~4.1]{BD25integral}, for $\N(\afr_1\afr_2)$ confined to a dyadic interval, the contribution is $O(B(\log B)^{1+2q})$, which is satisfactory after summing over the $O(\log T)$ many dyadic intervals.

    Similarly, when $\max_v |a_6|_v$ or $\max_v |a_7|_v$ is confined to a dyadic interval, their corresponding number is $O((\log B)^{q-1})$ by \cite[Lemma~3.1]{BD25integral} (in the extremal case $q=0$, note that their number is in fact $0$ for sufficiently large $B$).  The sum over $a_1,a_2$ (now without a restriction) yields another factor of $B(\log B)^2$, leaving us again with a satisfactory bound.
\end{proof}

\section{Counting II: The main contribution}\label{sec:main_contribution}

Having prepared the stage for the counting argument, we complete the proof of Theorem~\ref{thm:main_concrete} in this section. In Section~\ref{sec:o-minimality}, we evaluate the lattice point count for $a_8$ by describing everything through an o-minimal structure, so that the Lipschitz principle is applicable.

This leads us to a main term involving a volume which we transform into archimedean densities in Section~\ref{sec:archimedean_densities}, before summing the arithmetic part over the remaining variables $a_1,a_2$ and the units $a_6, a_7$ in Section~\ref{sec:endgame}, leading us to the desired asymptotic formula.

\subsection{Counting via o-minimal structures}\label{sec:o-minimality}

Recall the definition of $\Gs$ from Section~\ref{sec:dependent_variable}.

\begin{lemma}\label{lem:lattice_points}
    For $\ab' \in \Os'_*$ with $\theta_0(\afrb')=1$ and $a_6,a_7 \in \OK^\times$, we have
    \begin{align*}
        |\Gs(\ab',a_6,a_7,\cfrb) \cap \Fs_0(\ab',a_6,a_7;B)| ={}& \frac{2^{r_2}\vol S_F(\ab',(a_6^{(v)},a_7^{(v)})_v;B)}{|\Delta_K|^{1/2}\N(\afr_1\Os_8)}\\ &+ O\left(\frac{B}{T^{1/d}\N(\afr_1\afr_2)}\right),
    \end{align*}
\end{lemma}

\begin{proof}
    Let
    \begin{equation*}
        \tau=(\tau_v)_v : \prod_{v \mid \infty} K_v \to \prod_{v \mid \infty} K_v,\quad (x_{8v})_v \mapsto (x_{8v}/\N(\afr_1\Os_8)^{1/d})_v.
    \end{equation*}
    Defining $S_F' = S_F-\sigma(\gamma)$ and $\Lambda = \tau(\sigma(\Gs'))$, we have $|\Gs\cap \Fs_0| = |\Lambda \cap \tau(S_F')|$. By construction, $\Lambda$ is a lattice of rank $d$ and determinant $2^{-r_2}|\Delta_K|^{1/2}$ with first successive minimum $\lambda_1\ge 1$. As in \cite[Lemma~9.1]{FP16} (modified as in \cite[Lemma~4.4]{BD25integral} compared to \cite[Lemma~5.5]{BD24}), $\tau(S_F')$ is the fiber of a set definable in the o-minimal structure $\RR_{\exp}$. Hence we can apply \cite[Theorem~1.3]{BW14}, which gives the main term as stated and an error term bounded (up to a bounded constant) by the sum over all proper coordinate subspaces $W$ of $\RR^d$ of the $\dim W$-dimensional volumes $V_W$ of the orthogonal projections of $\tau(S_F')$ to $W$. These volumes remain the same when replacing $\tau(S_F')$ by $\tau(S_F)$.

    Defining
    \begin{equation*}
        S_F^{(v)} = \{x_{8v} \in K_v : |x_{8v}|_v \le c_v |B_{8v}|_v\}
    \end{equation*} for sufficiently large $c_v>0$, we
    observe using \eqref{eq:bound_ai} that $\tau(S_F) \subset \prod_{v \mid \infty} \tau_v(S_F^{(v)})$. For $p_v \in \{0,\dots,d_v\}$, the orthogonal projection of $\tau_v(S_F^{(v)})$ to a $p_v$-dimensional coordinate subspace of $K_v \cong \RR^{d_v}$ has $p_v$-dimensional volume
    \begin{equation*}
        V^{(v)}_{p_v} \ll \left(\frac{B}{\N(\afr_1\afr_2)}\right)^{p_v/d}
    \end{equation*}
    since $\tau(S_F^{(v)})$ is contained in a box of side length $\ll (B/(\N(\afr_1\afr_2)))^{1/d}$.
    
    Hence $V_W \le \prod_{v \mid \infty} V^{(v)}_{p_v}$ for certain $p_v \in \{0,\dots,d_v\}$ that are not all zero. Therefore,
    \begin{equation*}
        V_W \ll \left(\frac{B}{\N(\afr_1\afr_2)}\right)^{1-1/d} \ll \frac{B}{T^{1/d}\N(\afr_1\afr_2)},
    \end{equation*}
    using \eqref{eq:restrict_a1a2} in the final step.
\end{proof}

\subsection{Archimedean densities}\label{sec:archimedean_densities}

Here, we compute the volume of the set $S_F$ appearing in the main term in Lemma~\ref{lem:lattice_points}. It grows linearly in $B$, and essentially the archimedean densities appear in the leading constant.

\begin{lemma}\label{lem:vol_S_F}
    For $\ab' \in \Os_*' \cap \Fs_1^5$ with $\theta_0(\afrb')=1$ and $a_6,a_7 \in \OK^\times$ with \eqref{eq:restrict_a6_a7}, we have
    \begin{equation*}
        \vol S_F(\ab',(a_6^{(v)},a_7^{(v)})_v;B) = \frac{2^{r_1}\pi^{r_2}R_K B}{|N(a_2)|} \left(1 + O(T^{-1/d})\right).
    \end{equation*}
\end{lemma}

\begin{proof}
    We substitute $x_{8v} = B_{8v}z_{8v}$ (with Jacobian $|\prod_{v \mid \infty} B_{8v}| = B/|N(a_2)|$) in the definition of $S_F$, and introduce the notation $a_i^{(v)} = B_{iv} z_{iv}$ for $i=6,7$. Then 
    \begin{equation*}
        \Nt_v(\ab',a_6^{(v)},a_7^{(v)},x_{8v})/B=\Nt_v(\oneb,z_{6v},z_{7v},z_{8v}),
    \end{equation*}
    and hence
    \begin{equation*}
        \vol S_F(\ab',(a_6^{(v)},a_7^{(v)})_v;B) = \frac{B}{|N(a_2)|} \vol S_F(\oneb, (z_{6v},z_{7v})_v; 1).
    \end{equation*}
    We have
    \begin{align*}
        &\vol S_F(\oneb, (z_{6v},z_{7v})_v; 1) 
        = \int_{(\Nt_v(\oneb,z_{6v},z_{7v},z_{8v}))_{v \mid \infty} \in \exp(F(1))} \prod_{v\mid \infty} \ddd z_{8v}\\
        &= \int_{(\Nt_v(\oneb,0,0,z_{8v}))_{v \mid \infty} \in \exp(F(1))} \prod_{v\mid \infty} \ddd z_{8v}
        + O\left(\int_{(\Nt_v(\oneb,0,0,z_{8v}))_{v \mid \infty} \in D} \prod_{v\mid \infty} \ddd z_{8v}\right)
    \end{align*}
    where $D$ is the $\delta$-neighborhood of the boundary of $\exp(F(1))$ for $\delta \ll T^{-1/d}$ (since $\Nt_v(\oneb,z_{6v},z_{7v},z_{8v}) = \Nt_v(\oneb,0,0,z_{8v})+O(T^{-d_v/d})$ by \eqref{eq:restrict_a6_a7}). 
    
    We complete the proof following \cite[Lemma~5.1]{FP16}. We define $f : \prod_{v \mid \infty} K_v \to \RR_{\ge 0}^{\Omega_\infty}$ as $(z_{8v})_v \mapsto (\Nt_v(\oneb,0,0,z_{8v}))_v$. Then the main term is $f_*(\vol)(\exp(F(1)))$, with
    \begin{equation*}
        f_*(\vol) = \prod_{v \mid \infty} \vol\{ z_{8v} \in K_v : \Nt_v(\oneb,0,0,z_{8v}) \le 1\} \cdot \vol = 2^{r_1} \pi^{r_2} \cdot \vol
    \end{equation*}
    (since this volume is simply the volume of a neighborhood of the origin with radius $1$) and $\vol(\exp(F(1))) = R_K$. Furthermore, the error term is $O(f_*(\vol)(D)) = O(\vol(D)) = O(T^{-1/d})$, using \cite[Lemma~4.6]{BD25integral}.
\end{proof}

\subsection{Completion of the proof}
\label{sec:endgame}

In this section, we complete the proof by performing the summations over $a_1,a_2$ (rewritten as summations over ideals) and over the units $a_6,a_7$, which result in the $\pfr$-adic densities and the factors $\alpha$ and $(\log B)^{2+2q}$.

\begin{prop}\label{prop:elements_to_ideals}
    We have
    \begin{equation*}
        \begin{aligned}
            \sum_{\cfrb \in \Cs} |M_\cfrb(B)| ={}& \frac{2^{r_1}(2\pi)^{r_2}(1 + O(T^{-1/d})) R_K}{|\Delta_K|^{1/2}} \sums{(\afr_1,\afr_2) \in \IK^2\\\eqref{eq:restrict_a1a2}\\        [\afr_1\afr_2]=[\OK]} \frac{\theta_0(\afrb')B}{\N(\afr_1\afr_2)} \sums{a_6,a_7 \in \OK^\times\\\eqref{eq:restrict_a6_a7}} 1\\
        &+ O(B(\log B)^{1+2q}(\log\log B))
        \end{aligned}
    \end{equation*}
\end{prop}

\begin{proof}
    We combine Proposition~\ref{prop:restrict_a6_a7} with \eqref{eq:lattice_points}, Lemma~\ref{lem:lattice_points}, and Lemma~\ref{lem:vol_S_F}. Recall that $\ab' \in \Os_*'\cap \Fs_1^5$ implies $a_3=a_4=a_5=1$.
    
    Summing the error term from Lemma~\ref{lem:lattice_points} over the remaining variables $a_1,\dots,a_7$ yields a total contribution bounded by
    \begin{equation*}
        \ll \frac{1}{T^{1/d}} \sums{\ab' \in \Os'_* \cap \Fs_1^5\\\eqref{eq:restrict_a1a2}} \sums{a_6,a_7 \in \OK^\times\\ \eqref{eq:restrict_a6_a7}} \frac{B}{\N(\afr_1\afr_2)} \ll \frac{B(\log B)^{2+2q}}{T^{1/d}}
    \end{equation*}
    by an argument similar to the proof of Proposition~\ref{prop:restrict_a6_a7}, which is satisfactory.
   
    In the main term, we replace the summations over $(a_1,a_2) \in (\Os_{1*}\times\Os_{2*})\cap \Fs_1^2$ and $\cfrb \in \Cs$ by a summation over ideals $\afr_1,\afr_2$ such that $\afr_1\afr_2$ is a principal ideal (since $\Os_1\Os_2 = \OK$), and we use that $|N(a_2)|\N(\Os_8) = \N(\afr_2)$ (since $\Os_8 = \Os_2^{-1}$).
\end{proof}

\begin{lemma}\label{lem:sum_a1_a2}
    We have
    \begin{align*}
        \sums{(\afr_1,\afr_2) \in \IK^2\\\eqref{eq:restrict_a1a2} \\      [\afr_1\afr_2]=[\OK]} \frac{\theta_0(\afrb')B}{\N(\afr_1\afr_2)} \sums{a_6,a_7 \in \OK^\times\\\eqref{eq:restrict_a6_a7}} 1 ={}& \frac{|\mu_K|^2}{h_K R_K^2} \rho_K^2 \alpha \theta_1 B(\log B)^{2+2q}\\&+O(B(\log B)^{1+2q}(\log\log B))
    \end{align*}
    with
    \begin{equation*}
        \theta_1 = \prod_\pfr \left(1-\frac{1}{\N\pfr^2}\right).
    \end{equation*}
\end{lemma}

\begin{proof}
    By \cite[Lemma~3.1]{BD25integral}, the sum over $a_6,a_7$ is
    \begin{equation*}
        \frac{|\mu_K|^2}{R_K^2 q!^2}(\log(B_6/T)\log(B_7/T))^q + O((\log B)^{2q-1}).
    \end{equation*}
    Since $|\theta_0(\afrb')| \le 1$, the sum of the error term over $\afr_1,\afr_2$ is $O(B(\log B)^{1+2q})$.

    For the main term, because of \eqref{eq:B6B7B8} and \eqref{eq:restrict_a1a2}, we consider
    \begin{equation*}
        V(z_1,z_2;B)\coloneq\frac{B(\log(B/T)-\log(z_1z_2))^q(\log(B/T))^q}{(\log B)^{2q}z_1z_2}\cdot V'(z_1,z_2;B),
    \end{equation*}
    where $V'$ is the indicator function of the set of $z_1,z_2\ge 1$ with $z_1z_2 \le B/T$.

    Since $V \ll B/(z_1z_2)$ vanishes unless $z_1,z_2 \ll B$, we can apply a version of \cite[Proposition~7.2]{DF14} where, for fixed $\afr_1$, the summation over $\afr_2$ is restricted to the ideal class inverse to that of $\afr_1$, which leads to an extra factor of $h_K^{-1}$. (This modification of \cite[Proposition~7.2]{DF14} is straightforward: \cite[Lemma~2.5]{DF14} provides a formula for the average growth of certain arithmetic functions when summed over a fixed ideal class, hence such a restriction is easily incorporated into \cite[Corollary~2.7]{DF14}, which in turn is used in the computation of the main term in \cite[Proposition~7.2]{DF14}.) Defining 
    \begin{equation*}
        V_0(B) = \frac{1}{q!^2} \ints{z_1,z_2\ge 1\\z_1z_2 \le B/T} \frac{B(\log(B/T)-\log(z_1z_2))^q(\log(B/T))^q}{z_1z_2} \ddd z_1 \ddd z_2,
    \end{equation*}
    we obtain the main term
    \begin{equation*}
        \frac{|\mu_K|^2 \rho_K^2 \theta_1}{h_K R_K^2} V_0(B) + O(B(\log B)^{1+2q}(\log \log B))
    \end{equation*}
    since the average as in \cite[Lemma~2.8]{DF14} of $\theta_0(\afrb')$ over $\afr_1,\afr_2$ is $\theta_1$ as in our statement.

    Substituting $z_i=(B/T)^{t_i}$, the range of integration turns into $\{(t_1,t_2) \in \RR_{\ge 0}^2 : t_1+t_2 \le 1\}$, and
    \begin{equation*}
        \log(B/T)-\log(z_1z_2)=\log(B/T)(1-t_1-t_2).
    \end{equation*}
    Hence, by definition of $\alpha$ and $T$, we obtain
    \begin{equation*}
        V_0(B) = \alpha B(\log(B/T))^{2+2q} = \alpha B(\log B)^{2+2q}+O(B(\log B)^{1+2q}(\log \log B)).\qedhere
    \end{equation*}
\end{proof}

Finally, the main result is obtained by combining Proposition~\ref{prop:parameterization}, Proposition~\ref{prop:elements_to_ideals}, Lemma~\ref{lem:sum_a1_a2}, and \eqref{eq:def_rho_K}.

\bibliographystyle{alpha}

\bibliography{dp4}

\end{document}